# TWO GROUPS $2^3.PSL_2(7)$ AND $2^3:PSL_2(7)$ OF ORDER 1344

Mehmet Koca[1], Ramazan Koç[2] and Nazife Ozdes Koca[3,*]

[2]*Department of Physics, Gaziantep University, 27310, Gaziantep, Turkey*

[3]*Department of Physics, College of Science, Sultan Qaboos University, Muscat, Oman*

*Correspondence e-mail: nazife@squ.edu.om

**Abstract**

We analyze the group structures of two groups of order 1344 which are respectively non-split and split extensions of the elementary Abelian group of order 8 by its automorphism group $PSL_2(7)$. They share the same character table. The group $2^3.PSL_2(7)$ is a finite subgroup of the Lie Group $G_2$ preserving the set of octonions $\pm e_i, (i = 1,2,...,7)$ representing a 7-dimensional octahedron. Its three maximal subgroups $2^3:7:3$, $2^3.S_4$ and $4.S_4:2$ correspond to the finite subgroups of the Lie groups $G_2, SO(4)$ and $SU(3)$ respectively. The group $2^3:PSL_2(7)$ representing the split extension possesses five maximal subgroups $2^3:7:3$, $2^3:S_4$, $4:S_4:2$ and two non-conjugate Klein's group $PSL_2(7)$. The character tables of the groups and their maximal subgroups, tensor products and decompositions of the irreducible representations under the relevant maximal subgroups are identified. Possible implications in physics are discussed.

## I. INTRODUCTION

We have introduced some part of this work in an earlier publication [1]. Since then we have observed that the simple group like $PSL_2(7)$ [2] and the subgroups thereof $7:3$ [3], $S_4$ and $A_4$ [4, 5, 6, 7, 8, 9] have been proposed to explain the properties of the Tri-Bimaximal Neutrino Mixing [10]. When the charged leptons and quark masses are incorporated into the scheme we may think of much larger discrete symmetries broken to the aforementioned finite subgroups of $SU(3)$. Along with this line we would like to introduce two groups of order 1344 which are non-split and split extensions of the elementary Abelian group $2^3$ of order 8 by its automorphism group $PSL_2(7)$.

[1] mehmetkocaphysics@gmail.com
[2] koc@gantep.edu.tr
[3] nazife@squ.edu.om

Before we proceed further a glossary may be introduced for the group theoretical concepts and notations used throughout the paper. We follow the notations of the Atlas of Finite Groups [11]. A cyclic group of order $p$ is denoted by $p$. An elementary Abelian group of order $p^n$ (denoted also by $p^n$) is the direct product of $n$ cyclic groups of each having order $p$. Thus, the elementary abelian group $2^3$ is the direct product of three cyclic groups of order $2$. The group $A.B$ denotes any group possessing a normal subgroup $A$, for which the corresponding quotient group is $B$. This is used for the non-split extension. The group $A{:}B$ indicates the split extension, or the semi-direct product $A \rtimes B$. Here a copy of the quotient group $B$ is a subgroup of the group $A{:}B$. This shows that the intersection of $A$ and $B$ is just the unit element. The quotient group of interest here is the special projective group $PSL_2(7)$ and has been discussed in the physics literature extensively [2,12,13].

The paper is organized as follows. For the group $2^3.PSL_2(7)$ is the automorphism group of the octonionic set $\pm e_i, (i = 1,2, ... ,7)$ we review in Section 2 the basic structure of the octonion algebra and introduce the group $7{:}3$ as the automorphism of the 7 imaginary units $e_i, (i = 1,2, ... ,7)$. In Section 3 we extend the automorphism to the full group of automorphism including the change of sign of the imaginary units of octonions and point out that the diagonal matrices form the elementary abelian subgroup of the automorphism group. The quotient group $PSL_2(7)$ is explicitly constructed and the maximal subgroups of the group $2^3.PSL_2(7)$ are identified. In Section 4 we discuss the construction of the group by using a finite subgroup of $SO(4)$ based on the preservation of the quaternion subalgebra of the octonion algebra. Section 5 deals with the construction of the 7-dimensional irreducible representation of the group $2^3{:}PSL_2(7)$ and its five maximal subgroups, two of which, are the non-conjugate Klein's group $PSL_2(7)$. In concluding Section 6 we point out as to how these groups can be used in physics. In Appendix **A** we study the tensor products of the irreducible representations. Appendix **B** lists the decompositions of the irreducible representations under the maximal subgroups.

## II. OCTONIONS AND THE GROUP $7{:}3$

The octonions (Cayley numbers) are sets of real numbers

$$q = (q_0, q_1, q_2, q_3, q_4, q_5, q_6, q_7)$$
$$= q_0.1 + q_1 e_1 + q_2 e_2 + q_3 e_3 + q_4 e_4 + q_5 e_5 + q_6 e_6 + q_7 e_7 \qquad (1)$$

where $e_i, (i = 1,2, .... ,7)$ are 7 imaginary octonionic units. Octonions are added like vectors and multiplied as follows

$$1. e_i = e_i.1 = e_i,$$

$$e_i e_j = -\delta_{ij} + \sum_{k=1}^{7} \emptyset_{ijk}\, e_k \quad (i,j,k = 1,2, ... ,7) \qquad (2)$$

where $\emptyset_{ijk}$ are completely antisymmetric in $i,j,k$ with the values $\pm 1$. We choose the basis as shown in Figure 1 such that [14]

$$\emptyset_{123} = \emptyset_{246} = \emptyset_{435} = \emptyset_{367} = \emptyset_{651} = \emptyset_{572} = \emptyset_{714} = 1, \qquad (3)$$

which follows from the cyclic rotation of the triangle in Figure 1.

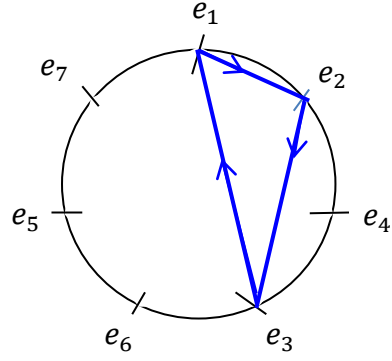

Figure 1. Octonionic multiplication based on quaternionic multiplication

The 7-imaginary units form 35 triads, 7 of which are associative and follow the ordering such that $e_i e_j e_k = -1$ when $i, j, k$ take one of the values in (3).
The 28 anti-associative triads can be obtained from the following four anti-associative triads

$$\begin{aligned} e_1(e_2 e_4) &= e_5 = -(e_1 e_2)e_4 \\ -e_1(e_2 e_6) &= e_7 = (e_1 e_2)e_6 \\ -e_1(e_2 e_5) &= e_4 = (e_1 e_2)e_5 \\ e_1(e_2 e_7) &= e_6 = -(e_1 e_2)e_7 \end{aligned} \qquad (4)$$

by using the cyclic permutation (1243657). The seven associative triads correspond to the lines of the finite projective geometry of seven lines and seven points of the Fano plane.

Without changing the signs of the octonionic imaginary units in Figure 1 we can obtain two operations preserving the octonion algebra. Let $\alpha$ transform the octonionic units $e_i$ in the cyclic order $(e_1 e_2 e_4 e_3 e_6 e_5 e_7)$ so that $\alpha^7 = 1$. Let $\beta$ fix $e_1$ and permutes the associative triad $(e_2 e_4 e_6)$ and the anti-associative triad $(e_3 e_7 e_5)$ in the indicated orders. It is clear that $\beta^3 = 1$ and preserves the octonion algebra. Indeed, the transformation $(e_3 e_7 e_5)$ is sufficient to determine the rest of the transformations of the quaternionic imaginary units. It is easy to prove that two generators generate a finite group of order 21 satisfying the generation relations

$$\alpha^7 = \beta^3 = \beta^{-1} \alpha \beta \alpha^3 = 1. \qquad (5)$$

It is clear from (5) that $\alpha$ generates an invariant subgroup of order 7 and hence the structure of the group $7:3$. The conjugacy classes are given by the set of group elements

$$\begin{aligned} C_1 &= (1), C_2 = (\alpha, \alpha^2, \alpha^4), C_3 = (\alpha^3, \alpha^5, \alpha^6), \\ C_4 &= (\alpha^a \beta), C_5 = (\alpha^a \beta^2), (a = 0, 1, \ldots, 6). \end{aligned} \qquad (6)$$

The character table of the group is depicted in TABLE I.

TABLE I. The character table of the Frobenius group $7:3$

| $7:3$ | $C_1$ | $7C_2^{[3]}$ | $7C_3^{[3]}$ | $3C_4^{[7]}$ | $3C_5^{[7]}$ |
|---|---|---|---|---|---|
| $\chi^{[1]}$ | 1 | 1 | 1 | 1 | 1 |
| $\chi^{[1_1]}$ | 1 | $\mu$ | $\bar{\mu}$ | 1 | 1 |
| $\chi^{[1_2]}$ | 1 | $\bar{\mu}$ | $\mu$ | 1 | 1 |
| $\chi^{[3_1]}$ | 3 | 0 | 0 | $\eta$ | $\bar{\eta}$ |
| $\chi^{[3_2]}$ | 3 | 0 | 0 | $\bar{\eta}$ | $\eta$ |

Here, $\chi^{[1_2]} = \chi^{[\bar{1_1}]}$, $\chi^{[3_2]} = \chi^{[\bar{3_1}]}$, $\mu = \frac{1}{2}(-1 + i\sqrt{3})$, $\eta = \frac{1}{2}(-1 + i\sqrt{7})$.

One can easily check that the 7-dimensional representation obtained by the generators $\alpha$ and $\beta$ is reducible and can be decomposed as $7 = 1 + 3_1 + \overline{3_1}$.

## III. THE GROUP $2^3.PSL_2(7)$ AS THE AUTOMORPHISM GROUP OF THE OCTONIONIC SET $\pm e_i$

Now we allow the octonionic units $e_i$ take also negative values, in other words, we also include the octonionic conjugates in the transformations. Any linear transformation on a non-associative triad determines the transformations of 7 imaginary units. Let us take the non-associative triad $(e_1 e_2 e_7)$ and impose the transformation $e_1 \to e_1$, $e_2 \to e_2$, $e_7 \to -e_7$. This transformation leads to the transformations $e_3 \to e_3, e_4 \to -e_4, e_5 \to -e_5$ and $e_6 \to -e_6$. Let us call it
$N_1 = (e_1 \to e_1, e_2 \to e_2, e_3 \to e_3, e_4 \to -e_4, e_5 \to -e_5, e_6 \to -e_6, e_7 \to -e_7)$ and denote it as a diagonal matrix with non-zero entities

$$N_1 = (1, 1, 1, -1, -1, -1, -1). \tag{7}$$

This diagonal transformation leaves the associative triad (123) intact and changes the signs of the other 4 octonionic units. When we take $\alpha$, $\beta$ and $N_1$ as generators we obtain a group of order 168 but not isomorphic to $PSL_2(7)$. It has, by construction, an invariant subgroup of order 8 generated by three diagonal transformations $N_1, N_2$ and $N_7$ where $N_2 = (1, -1, -1, 1, -1, -1, 1)$ and $N_7 = (-1, 1, -1, 1, -1, 1, -1)$. They generate the elementary Abelian group $2^3$ where the other elements are defined as

$$\begin{aligned} N_3 &= N_1 N_2 = N_2 N_1, \quad N_4 = N_7 N_1 = N_1 N_7, \\ N_5 &= N_7 N_2 = N_2 N_7, \quad N_6 = N_7 N_3 = N_3 N_7. \end{aligned} \tag{8}$$

More compactly, they can be written as

$$N_i N_j = N_j N_i = N_k, (ijk = 123, 147, 165, 246, 257, 345, 367). \tag{9}$$

The 7-diagonal matrices can be used to define the Fano plane. One can show that the set of group elements $N_i, (i = 1, 2, ..., 7)$ is invariant under the conjugations of the generators so that the group can be designated as $2^3: 7: 3$, the order of which is 168 but since it has an invariant subgroup it is not isomorphic to the simple group $PSL_2(7)$. The 7-dimensional representation of the group $2^3: 7: 3$ obtained from the generators $\alpha$, $\beta$ and $N_1$ is irreducible and denoted by $7_1$ in the character table of the group $2^3:7:3$ given in TABLE II.

TABLE II. The character table of the group $2^3:7:3$

| $2^3:7:3$ | $C_1$ | $7C_2^{[2]}$ | $28C_3^{[3]}$ | $28C_4^{[3]}$ | $28C_5^{[6]}$ | $28C_6^{[6]}$ | $24C_7^{[7]}$ | $24C_8^{[7]}$ |
|---|---|---|---|---|---|---|---|---|
| $\chi^{[1]}$ | 1 | 1 | 1 | 1 | 1 | 1 | 1 | 1 |
| $\chi^{[1_1]}$ | 1 | 1 | $\bar{\mu}$ | $\mu$ | $\bar{\mu}$ | $\mu$ | 1 | 1 |
| $\chi^{[1_2]}$ | 1 | 1 | $\mu$ | $\bar{\mu}$ | $\mu$ | $\bar{\mu}$ | 1 | 1 |
| $\chi^{[3_1]}$ | 3 | 3 | 0 | 0 | 0 | 0 | $\eta$ | $\bar{\eta}$ |
| $\chi^{[3_2]}$ | 3 | 3 | 0 | 0 | 0 | 0 | $\bar{\eta}$ | $\eta$ |
| $\chi^{[7_1]}$ | 7 | -1 | 1 | 1 | -1 | -1 | 0 | 0 |
| $\chi^{[7_2]}$ | 7 | -1 | $\mu$ | $\bar{\mu}$ | $-\mu$ | $-\bar{\mu}$ | 0 | 0 |
| $\chi^{[7_3]}$ | 7 | -1 | $\bar{\mu}$ | $\mu$ | $-\bar{\mu}$ | $-\mu$ | 0 | 0 |

Note that in TABLE II the characters satisfy the relations $\chi^{[1_2]} = \chi^{[\overline{1_1}]}, \chi^{[3_2]} = \chi^{[\overline{3_1}]}, \chi^{[7_3]} = \chi^{[\overline{7_1}]}$ and the other irreducible representations are real. The group $2^3:7:3$ is not the full automorphism group of the set of octonions $\pm e_i$. In fact, a transformation of the form $\gamma:(e_1 \leftrightarrow -e_4, e_2 \leftrightarrow -e_5, e_3 \to -e_3, e_6 \to e_6, e_7 \to -e_7)$ preserves the octonion algebra and does not belong to the group of elements of the group $2^3:7:3$. Adjoining $\gamma$ as a new generator then $\alpha, \beta, \gamma$ and $N_1$ together generate the group of order 1344. We will prove that it has the structure $2^3.PSL_2(7)$ which can also be generated by two generators $\alpha$ and $\gamma$ only. Here the group $PSL_2(7)$ is the quotient group whose generators are defined by the conjugations over the elements of the elementary abelian group $2^3$:

$$\tilde{\alpha}: N_i \to \alpha^{-1}N_i\alpha, \ \tilde{\beta}: N_i \to \beta^{-1}N_i\beta, \ \tilde{\gamma}: N_i \to \gamma^{-1}N_i\gamma. \qquad (10)$$

They permute the diagonal matrices as

$$\begin{aligned} \tilde{\alpha} &= (N_1 N_2 N_4 N_3 N_6 N_5 N_7), \\ \tilde{\beta} &= (N_3 N_2 N_1)(N_4 N_6 N_5)(N_7) \\ \tilde{\gamma} &= (N_1 N_5)(N_2)(N_3 N_7)(N_4)(N_6). \end{aligned} \qquad (11)$$

The group generated by the generators in (11) is isomorphic to the Klein's group $PSL_2(7)$ of order 168 and the 7-dimensional representation obtained from (11) is reducible $7 = 1 + 6$. This is expected because 7 diagonal matrices $N_i$ form the Fano plane whose automorphism group is the Klein's group. The character table of the group $PSL_2(7)$ is shown in TABLE III. (Note the difference between TABLE II and TABLE III).

TABLE III. The character table of the group $PSL_2(7)$

| $PSL_2(7)$ | $C_1$ | $21C_2^{[2]}$ | $56C_3^{[3]}$ | $42C_4^{[4]}$ | $42C_5^{[7]}$ | $42C_6^{[7]}$ |
|---|---|---|---|---|---|---|
| $\chi^{[1]}$ | 1 | 1 | 1 | 1 | 1 | 1 |
| $\chi^{[3_1]}$ | 3 | -1 | 0 | 1 | $\eta$ | $\bar{\eta}$ |
| $\chi^{[3_2]}$ | 3 | -1 | 0 | 1 | $\bar{\eta}$ | $\eta$ |
| $\chi^{[6]}$ | 6 | 2 | 0 | 0 | -1 | -1 |
| $\chi^{[7]}$ | 7 | -1 | 1 | -1 | 0 | 0 |
| $\chi^{[8]}$ | 8 | 0 | -1 | 0 | 1 | 1 |

Here again the characters satisfy $\chi^{[3_2]} = \chi^{\overline{[3_1]}}$. A copy of the $PSL_2(7)$ does not exist in the group $2^3.PSL_2(7)$, therefore it is not a subgroup and hence the latter group is the non-split extension of the group $2^3$ by the group $PSL_2(7)$. The character table of the groups $2^3.PSL_2(7)$ and $2^3:PSL_2(7)$ is given in TABLE IV.

TABLE IV. The character table of the groups $2^3.PSL_2(7)$ and $2^3:PSL_2(7)$

| $2^3.PSL_2(7)$ | $C_1$ | $C_2^{[2]}$ | $C_3^{[2]}$ | $C_4^{[3]}$ | $C_5^{[4]}$ | $C_6^{[4]}$ | $C_7^{[6]}$ | $C_8^{[7]}$ | $C_9^{[7]}$ | $C_{10}^{[8]}$ | $C_{11}^{[8]}$ |
|---|---|---|---|---|---|---|---|---|---|---|---|
| $2^3:PSL_2(7)$ | $C_1$ | $C_2^{[2]}$ | $C_3^{[4]}$ | $C_4^{[3]}$ | $C_5^{[2]}$ | $C_6^{[2]}$ | $C_7^{[6]}$ | $C_8^{[7]}$ | $C_9^{[7]}$ | $C_{10}^{[4]}$ | $C_{11}^{[4]}$ |
| # of elements | 1 | 7 | 84 | 224 | 42 | 42 | 224 | 192 | 192 | 168 | 168 |
| $\chi^{[1]}$ | 1 | 1 | 1 | 1 | 1 | 1 | 1 | 1 | 1 | 1 | 1 |
| $\chi^{[3_1]}$ | 3 | 3 | -1 | 0 | -1 | -1 | 0 | $\eta$ | $\bar{\eta}$ | 1 | 1 |
| $\chi^{[3_2]}$ | 3 | 3 | -1 | 0 | -1 | -1 | 0 | $\bar{\eta}$ | $\eta$ | 1 | 1 |
| $\chi^{[6]}$ | 6 | 6 | 2 | 0 | 2 | 2 | 0 | -1 | -1 | 0 | 0 |
| $\chi^{[7_1]}$ | 7 | -1 | -1 | 1 | -1 | 3 | -1 | 0 | 0 | 1 | -1 |
| $\chi^{[7_2]}$ | 7 | 7 | -1 | 1 | -1 | -1 | 1 | 0 | 0 | -1 | -1 |
| $\chi^{[7_3]}$ | 7 | -1 | -1 | 1 | 3 | -1 | -1 | 0 | 0 | -1 | 1 |
| $\chi^{[8]}$ | 8 | 8 | 0 | -1 | 0 | 0 | -1 | 1 | 1 | 0 | 0 |
| $\chi^{[14]}$ | 14 | -2 | -2 | -1 | 2 | 2 | 1 | 0 | 0 | 0 | 0 |
| $\chi^{[21_1]}$ | 21 | -3 | 1 | 0 | -3 | 1 | 0 | 0 | 0 | -1 | 1 |
| $\chi^{[21_2]}$ | 21 | -3 | 1 | 0 | 1 | -3 | 0 | 0 | 0 | 1 | -1 |

The character table shows that all the representations are real except the 3-dimensional representations which satisfy $\chi^{[3_2]} = \chi^{\overline{[3_1]}}$.

The group $2^3.PSL_2(7)$ has three maximal subgroups having the structures $2^3:7:3$, $4.S_4:2$ and $2^3.S_4$ of respective orders 168, 192 and 192. The group $2^3:7:3$ preserves the octonion multiplication and the 7-dimensional representation denoted by $7_1$ is irreducible. It is a subgroup of the Lie group $G_2$ and the groups $4.S_4:2$ and $2^3.S_4$ are finite subgroups of the groups $SU(3)$ and $SO(4)$ respectively as we will study in the next section. The Lie groups $SU(3)$ and $SO(4)$ are the maximal subgroups of the Lie group $G_2$.

## Maximal subgroups $4.S_4:2$ and $2^3.S_4$ of the group $2^3.PSL_2(7)$

The group $2^3.PSL_2(7)$ has two maximal subgroups of order 192. Now we discuss some properties of the first subgroup $4.S_4:2$. Let $\bar{e}_i = -e_i$ denote the octonionic conjugate. Define the transformation

$$\theta = (e_1\bar{e}_5)(e_2\bar{e}_3 e_4\bar{e}_7\bar{e}_2 e_3\bar{e}_4 e_7)(e_6\bar{e}_6), \quad \text{with } \theta^8=1. \tag{12}$$

The generators $\gamma$ and $\theta$ generate the group $4.S_4:2$ which represents a finite subgroup of $SU(3)$ as the generators fix the octonion $\pm e_6$. The character table of the group $4.S_4:2$ is shown in TABLE V. This is also one of the maximal subgroups of the Chevalley's group $G_2(2)$[15] preserving the octonionic root system of $E_7$. The 7-dimensional representation can be written as $7 = 1 + 6$; more properly, $7_1 = 1_1 + 6_3$. By adjoining the generator $\alpha$ which permutes the 7 imaginary octonionic units one can generate the group $2^3.PSL_2(7)$. It is easy to prove that $\tilde{\gamma}$ and $\tilde{\theta} = (N_1 N_4 N_2 N_5)(N_3)(N_6 N_7)$ with $\tilde{\gamma}^2 = \tilde{\theta}^4 = 1$ generate the subgroup $S_4$ as the quotient group. Let $a$ and $b$ be the generators of a group, then the standard generation relations of $S_4$ is $a^4 = b^3 = (ab)^2 = 1$ [16]. This relation can be satisfied if we define $a = \tilde{\gamma}\tilde{\theta}\tilde{\gamma}$ and $b = \tilde{\gamma}\,\tilde{\theta}^{-1}$.

TABLE V. The character table of the groups $4.S_4:2$ and $4:S_4:2$

| $4.S_4:2$ | $C_1$ | $C_2^{[2]}$ | $C_3^{[2]}$ | $C_4^{[2]}$ | $C_5^{[2]}$ | $C_6^{[2]}$ | $C_7^{[3]}$ | $C_8^{[4]}$ | $C_9^{[4]}$ | $C_{10}^{[4]}$ | $C_{11}^{[4]}$ | $C_{12}^{[6]}$ | $C_{13}^{[8]}$ | $C_{14}^{[8]}$ |
|---|---|---|---|---|---|---|---|---|---|---|---|---|---|---|
| $4:S_4:2$ | $C_1$ | $C_2^{[2]}$ | $C_3^{[2]}$ | $C_4^{[4]}$ | $C_5^{[4]}$ | $C_6^{[4]}$ | $C_7^{[3]}$ | $C_8^{[2]}$ | $C_9^{[2]}$ | $C_{10}^{[2]}$ | $C_{11}^{[2]}$ | $C_{12}^{[6]}$ | $C_{13}^{[4]}$ | $C_{14}^{[4]}$ |
| # of elements | 1 | 3 | 4 | 12 | 12 | 12 | 32 | 12 | 6 | 6 | 12 | 32 | 24 | 24 |
| $\chi^{[1]}$ | 1 | 1 | 1 | 1 | 1 | 1 | 1 | 1 | 1 | 1 | 1 | 1 | 1 | 1 |
| $\chi^{[1_1]}$ | 1 | 1 | -1 | 1 | -1 | -1 | 1 | 1 | 1 | 1 | -1 | -1 | 1 | -1 |
| $\chi^{[1_2]}$ | 1 | 1 | -1 | -1 | 1 | -1 | 1 | -1 | 1 | 1 | 1 | -1 | -1 | 1 |
| $\chi^{[1_3]}$ | 1 | 1 | 1 | -1 | -1 | 1 | 1 | -1 | 1 | 1 | -1 | 1 | -1 | -1 |
| $\chi^{[2_1]}$ | 2 | 2 | -2 | 0 | 0 | -2 | -1 | 0 | 2 | 2 | 0 | 1 | 0 | 0 |
| $\chi^{[2_2]}$ | 2 | 2 | 2 | 0 | 0 | 2 | -1 | 0 | 2 | 2 | 0 | -1 | 0 | 0 |
| $\chi^{[3_1]}$ | 3 | 3 | 3 | -1 | -1 | -1 | 0 | -1 | -1 | -1 | -1 | 0 | 1 | 1 |
| $\chi^{[3_2]}$ | 3 | 3 | 3 | 1 | 1 | -1 | 0 | 1 | -1 | -1 | 1 | 0 | -1 | -1 |
| $\chi^{[3_3]}$ | 3 | 3 | -3 | 1 | -1 | 1 | 0 | 1 | -1 | -1 | -1 | 0 | -1 | 1 |
| $\chi^{[3_4]}$ | 3 | 3 | -3 | -1 | 1 | 1 | 0 | -1 | -1 | -1 | 1 | 0 | 1 | -1 |
| $\chi^{[6_1]}$ | 6 | -2 | 0 | 0 | -2 | 0 | 0 | 0 | -2 | 2 | 2 | 0 | 0 | 0 |
| $\chi^{[6_2]}$ | 6 | -2 | 0 | 0 | 2 | 0 | 0 | 0 | -2 | 2 | -2 | 0 | 0 | 0 |
| $\chi^{[6_3]}$ | 6 | -2 | 0 | -2 | 0 | 0 | 0 | 2 | 2 | -2 | 0 | 0 | 0 | 0 |
| $\chi^{[6_4]}$ | 6 | -2 | 0 | 2 | 0 | 0 | 0 | -2 | 2 | -2 | 0 | 0 | 0 | 0 |

The group $2^3.S_4$ can be generated by the generators

$$A = (e_1\bar{e}_7\ e_3\bar{e}_1 e_7\bar{e}_3)(e_2\bar{e}_4\bar{e}_6\bar{e}_2 e_4 e_6)(\ e_5\bar{e}_5),$$

$$B = (e_1)\ (e_2\bar{e}_6\bar{e}_2 e_6)(\ e_3\bar{e}_5\bar{e}_3 e_5)(e_4)(e_7), \quad (13)$$

$$A^6 = B^4 = 1.$$

The elements $\tilde{A}$ and $\tilde{B}$ can be obtained as follows

$$\tilde{A} = (N_1 N_6 N_2)(N_3 N_5 N_4)(N_7),\ \tilde{B} = (N_1 N_3)(N_2)(N_4 N_6)(N_5)(N_7). \quad (14)$$

Now define $a = \tilde{B}\tilde{A} = (N_1 N_4)(N_2 N_3 N_5 N_6)(N_7)$ and $b = \tilde{A}^{-1}$. It is straightforward to show that they satisfy the $S_4$ generation relations $a^4 = b^3 = (ab)^2 = 1$. Hence, the quotient group is $S_4$ as claimed. The matrix representations of $\tilde{A}$ and $\tilde{B}$ do not preserve the octonion algebra therefore a copy of $S_4$ does not exist in the group $2^3.S_4$ hence it is the non-split extension of the group $2^3$ by $S_4$.

It is clear from the generators $A$ and $B$ that they preserve the sets of octonionic units $(\pm e_2 \pm e_4 \pm e_6)$ and $(\pm e_1 \pm e_3 \pm e_5 \pm e_7)$ separately so that the $7 \times 7$ matrix representation can be written in the block diagonal form of $3 \times 3$ and $4 \times 4$ matrices. This proves that the 7-dimensional representation of the group $2^3.PSL_2(7)$ can be decomposed as $7_1 = 3_1 + 4_1$ under the group $2^3.S_4$.

A rigorous proof can also be given by noting that the group $2^3.S_4$ is a finite subgroup of $SO(4)$ which preserves the quaternion substructure of the octonion algebra.

Let $h$ represents a quaternion. Then an octonion can be written as $h + e_7 h$ where the quaternionic imaginary units are taken as $e_1, e_2$ and $e_3$. It can be proved that the quaternion preserving transformation on the octonion $h + e_7 h$ [17] can be written as

$$h + e_7 h \to ph\bar{p} + phq \quad (15)$$

where $p$ and $q$ are unit quaternions. The first term $ph\bar{p}$ preserves the scalar part of the octonion. Therefore, without loss of generality, we can also write it as $p\ Im(h\bar{p})$. Since we are interested in the transformations of the set $\pm e_i, (i = 1,2, ... ,7)$ we can write it as

$$Im\ (V_0) + e_7 V_0 \to pIm(\ V_0\bar{p}) + e_7 p V_0 q \quad (16)$$

where the set $V_0 = \{\pm 1, \pm e_1,\ \pm e_2, \pm e_3\}$ represents the elements of the quaternion group of order 8. So, the question now turns out to be that for which unit quaternions $p$ and $q$ the set of quaternions $V_0$ is preserved. The answer lies in the finite quaternion subgroups involving $V_0$ as a subgroup. They are the binary tetrahedral group $\mathcal{T}$, binary octahedral group $\mathcal{O}$ and the binary icosahedral group $\mathcal{J}$[16]. It turns out that we need to invoke the binary octahedral group whose sets of elements can be written as the union of 6 quaternionic sets [18]

$$\mathcal{O} = V_0 + V_+ + V_- + V_1 + V_2 + V_3. \quad (17)$$

Here the sets of quaternions are defined as

$$V_0 = \{\pm 1, \pm e_1, \pm e_2, \pm e_3\}$$
$$V_+ = \frac{1}{2}(\pm 1 \pm e_1 \pm e_2 \pm e_3),\ \text{even number of (+sign)}$$

$$V_- = \frac{1}{2}(\pm 1 \pm e_1 \pm e_2 \pm e_3), \text{ odd number of (+sign)}$$
$$V_1 = \{\frac{1}{\sqrt{2}}(\pm 1 \pm e_1), \frac{1}{\sqrt{2}}(\pm e_2 \pm e_3)\}, \tag{18}$$
$$V_2 = \{\frac{1}{\sqrt{2}}(\pm 1 \pm e_2), \frac{1}{\sqrt{2}}(\pm e_3 \pm e_1)\},$$
$$V_3 = \{\frac{1}{\sqrt{2}}(\pm 1 \pm e_3), \frac{1}{\sqrt{2}}(\pm e_1 \pm e_2)\}.$$

They satisfy the multiplication table as shown in TABLE VI.

TABLE VI. Multiplication table for the binary octahedral group

|       | $V_0$ | $V_+$ | $V_-$ | $V_1$ | $V_2$ | $V_3$ |
|-------|-------|-------|-------|-------|-------|-------|
| $V_0$ | $V_0$ | $V_+$ | $V_-$ | $V_1$ | $V_2$ | $V_3$ |
| $V_+$ | $V_+$ | $V_-$ | $V_0$ | $V_3$ | $V_1$ | $V_2$ |
| $V_-$ | $V_-$ | $V_0$ | $V_+$ | $V_2$ | $V_3$ | $V_1$ |
| $V_1$ | $V_1$ | $V_2$ | $V_3$ | $V_0$ | $V_+$ | $V_-$ |
| $V_2$ | $V_2$ | $V_3$ | $V_1$ | $V_-$ | $V_0$ | $V_+$ |
| $V_3$ | $V_3$ | $V_1$ | $V_2$ | $V_+$ | $V_-$ | $V_0$ |

To determine the pair of quaternions $p$ and $q$ we use the notation $phq \coloneqq [p,q]h$. Later, we simply drop the quaternion $h$ and write the $SO(4)$ group elements as the pair $[p,q]$. Since we have

$$V_0 = pV_0q \to q = \overline{V_0}\overline{p}V_0 = V_0\overline{p}V_0 \tag{19}$$

and $p$ can take values from the sets $V_0, V_+, V_-, V_1, V_2, V_3$ then the corresponding sets of quaternions $q$ can be determined from the TABLE VI as $V_0, V_-, V_+, V_1, V_2, V_3$. Therefore the $SO(4)$ group elements preserving the set $V_0$ can be written as the union of pairs of sets

$$[V_0,V_0] + [V_+,V_-] + [V_-,V_+] + [V_1,V_1] + [V_2,V_2] + [V_3,V_3]. \tag{20}$$

This is a group of order 192 with 13 conjugacy classes which can be converted to 4-dimensional irreducible representation if the basis is chosen as the unit quaternions $1, e_1, e_2, e_3$. When the term $[p, \bar{p}]$ is incorporated into the group $2^3.S_4$ it can be converted into $7 \times 7$ matrix representation in the block diagonal form of $3 \times 3$ and $4 \times 4$ matrices. This proves that the 7-dimensional irreducible representation of the group $2^3.PSL_2(7)$ can be decomposed as $7_1 = 3_1 + 4_1$ under the irreducible representations of the group $2^3.S_4$. The character table of the group $2^3.S_4$ is given in TABLE VII. Its elementary Abelian subgroup $2^3$ can be generated by the group elements, say, $[1, -1], [e_1, -e_1,]$ and $[e_2, -e_2,]$ and it can be proven that the group $2^3$ is a normal subgroup. It is also easy to identify its maximal subgroups. For example, the set $[V_0, V_0] + [V_+, V_-] + [V_-, V_+]$ forms a group of order 96 and $[V_0, V_0] + [V_1, V_1]$ forms a group of order 64.

TABLE VII. The character table of the groups $2^3.S_4$ and $2^3:S_4$

| $2^3.S_4$ | $C_1$ | $C_2^{[2]}$ | $C_3^{[2]}$ | $C_4^{[2]}$ | $C_5^{[2]}$ | $C_6^{[3]}$ | $C_7^{[4]}$ | $C_8^{[4]}$ | $C_9^{[4]}$ | $C_{10}^{[4]}$ | $C_{11}^{[6]}$ | $C_{12}^{[8]}$ | $C_{13}^{[8]}$ |
|---|---|---|---|---|---|---|---|---|---|---|---|---|---|
| $2^3:S_4$ | $C_1$ | $C_2^{[2]}$ | $C_3^{[2]}$ | $C_4^{[4]}$ | $C_5^{[4]}$ | $C_6^{[3]}$ | $C_7^{[2]}$ | $C_8^{[2]}$ | $C_9^{[2]}$ | $C_{10}^{[2]}$ | $C_{11}^{[6]}$ | $C_{12}^{[4]}$ | $C_{13}^{[4]}$ |
| # of elements | 1 | 1 | 6 | 12 | 24 | 32 | 6 | 6 | 12 | 12 | 32 | 24 | 24 |
| $\chi^{[1]}$ | 1 | 1 | 1 | 1 | 1 | 1 | 1 | 1 | 1 | 1 | 1 | 1 | 1 |
| $\chi^{[1_1]}$ | 1 | 1 | 1 | 1 | -1 | 1 | 1 | 1 | -1 | -1 | 1 | -1 | -1 |
| $\chi^{[2]}$ | 2 | 3 | 2 | 2 | 0 | -1 | 2 | 2 | 0 | 0 | -1 | 0 | 0 |
| $\chi^{[3_1]}$ | 3 | 3 | -1 | -1 | -1 | 0 | -1 | 3 | 1 | 1 | 0 | -1 | 1 |
| $\chi^{[3_2]}$ | 3 | 3 | -1 | -1 | 1 | 0 | -1 | 3 | -1 | -1 | 0 | 1 | -1 |
| $\chi^{[3_3]}$ | 3 | 3 | 3 | -1 | 1 | 0 | -1 | -1 | 1 | 1 | 0 | -1 | -1 |
| $\chi^{[3_4]}$ | 3 | 3 | -1 | -1 | -1 | 0 | 3 | -1 | 1 | 1 | 0 | 1 | -1 |
| $\chi^{[3_5]}$ | 3 | 3 | -1 | -1 | 1 | 0 | 3 | -1 | -1 | -1 | 0 | -1 | 1 |
| $\chi^{[3_6]}$ | 3 | 3 | 3 | -1 | -1 | 0 | -1 | -1 | -1 | -1 | 0 | 1 | 1 |
| $\chi^{[4_1]}$ | 4 | -4 | 0 | 0 | 0 | 1 | 0 | 0 | -2 | 2 | -1 | 0 | 0 |
| $\chi^{[4_2]}$ | 4 | -4 | 0 | 0 | 0 | 1 | 0 | 0 | 2 | -2 | -1 | 0 | 0 |
| $\chi^{[6]}$ | 6 | 6 | -2 | 2 | 0 | 0 | -2 | -2 | 0 | 0 | 0 | 0 | 0 |
| $\chi^{[8]}$ | 8 | -8 | 0 | 0 | 0 | -1 | 0 | 0 | 0 | 0 | 1 | 0 | 0 |

Since the group $2^3.S_4$ leaves the subsets $(\pm e_1, \pm e_2, \pm e_3)$ and $e_7(\pm 1, \pm e_1, \pm e_2, \pm e_3)$ invariant it is easy to construct its 7 conjugate groups by permuting the quaternionic subsets $(123), (246), (435), (367), (652), (572), (714)$ by the group generator $\alpha$ so that the group is extended by the inclusion of the generator $\alpha$ to the full automorphism group $2^3.PSL_2(7)$ of the octonionic set $\pm e_i$.

## IV. CONSTRUCTION OF THE 7-DIMENSIONAL IRREDUCIBLE REPRESENTATION OF THE GROUP $2^3:PSL_2(7)$

Although the groups are not isomorphic to each other, the split extension $2^3:PSL_2(7)$ shares the same character table with the non-split extension $2^3.PSL_2(7)$. A close inspection shows that the class structures are the same although the powers of the group elements are not always the same. To give an example from the lower rank groups possessing the same character table we may quote the dihedral group of order 8 and the quaternion group denoted by $V_0$ in our notation. The group $2^3:PSL_2(7)$ does not preserve the octonion algebra. Then let us assume that the group $2^3:PSL_2(7)$ acts on a 7-dimensional real vector space with the vector components $x_i (i = 1,2,...,7)$. Let us assume the notation $\overline{x_i} = -x_i$.

The construction of this group is easy by using its maximal subgroups. We have already discussed that the automorphism group $PSL_2(7)$ of the elementary abelian group $2^3$ can be generated by $\tilde{\alpha}, \tilde{\beta}$ and $\tilde{\gamma}$ in (11) and there exists a copy of the $PSL_2(7)$ in the group $2^3:PSL_2(7)$. For the representation of $\tilde{\alpha}, \tilde{\beta}$ and $\tilde{\gamma}$ defined over $N_i$, we can simply replace $N_i$ by $x_i$ and define the diagonal matrices $N_i$ as if they are acting on the vector components $x_i$. We can then simply adjoin an element from the elementary Abelian group $2^3$ to generate the whole group. Therefore, it would suffice to take $\widetilde{\alpha}, \widetilde{\beta}, \widetilde{\gamma}$ of $PSL_2(7)$ and $N_1$ as the generators of the group $2^3{:}PSL_2(7)$. Although this group does not preserve the octonion algebra there exists yet a maximal subgroup $2^3{:}7{:}3$ generated by $\widetilde{\alpha}, \widetilde{\beta}$ and $N_1$ preserving the octonion algebra as we discussed in Section 3. It is quite natural to expect that the groups $2^3{:}S_4$ and $4{:}S_4{:}2$ are maximal subgroups of order 192 in the group $2^3{:}PSL_2(7)$. The group $2^3{:}S_4$ can be generated by $\widetilde{A}, \widetilde{B}$ and $N_1$ and its character table is depicted in TABLE VII. Similarly, the subgroup $4{:}S_4{:}2$ can be generated by the generators $\widetilde{\gamma}$, $\widetilde{\theta}$ and $N_1$. Its character table is displayed in TABLE V.

The group $2^3{:}PSL_2(7)$ is a maximal subgroup of the simple group $A_8$, the even permutations of the 8 letters [11] isomorphic to the maximal rotation subgroup of the Coxeter-Weyl group $W(SU(8)) \cong S_8$.

We have already listed the maximal subgroups of the group $2^3{:}PSL_2(7)$ as the groups $PSL_2(7)$, $2^3{:}7{:}3$, $2^3{:}S_4$ and $4{:}S_4{:}2$. There exists yet another $PSL_2(7)$ not conjugate to the one generated by $\widetilde{\alpha}, \widetilde{\beta}$ and $\widetilde{\gamma}$ as Conway proved [19]. To see this, let us replace the generator $\widetilde{\gamma}$ by an another generator $\delta = (x_1\overline{x_5})(x_2)(x_3\overline{x_7})(x_4)(x_6)$, with $\delta \in 2^3{:}PSL_2(7)$ and $\delta^2 = 1$ instead of $\widetilde{\gamma}$. One can show that $\widetilde{\alpha}, \widetilde{\beta}, \delta$ and $N_1$ generate the same irreducible representation of the group $2^3{:}PSL_2(7)$ in which the $\widetilde{\alpha}, \widetilde{\beta}, \delta$ represent the generators of another group $PSL_2(7)$. This is a miraculous structure that the $PSL_2(7) = <\widetilde{\alpha}, \widetilde{\beta}, \delta>$ is not conjugate to the group $PSL_2(7) = <\widetilde{\alpha}, \widetilde{\beta}, \widetilde{\gamma}>$. It is straightforward to show that $\widetilde{\gamma}\delta = \delta\widetilde{\gamma} = N_7$. It is also true that 7-dimensional representation of the group $PSL_2(7) = <\widetilde{\alpha}, \widetilde{\beta}, \delta>$ is irreducible while the 7- dimensional representation of the group $PSL_2(7) = <\widetilde{\alpha}, \widetilde{\beta}, \widetilde{\gamma}>$ is reducible with $7 = 1 + 6$. Conway attributes this exceptional behavior to the holomorph of an elementary Abelian group to exhibit this miraculous feature.

Therefore, the group $2^3{:}PSL_2(7)$ has five maximal subgroups $2^3{:}7{:}3$, $2^3{:}S_4$, $4{:}S_4{:}2$ and two nonconjugate $PSL_2(7)$.

## V. CONCLUDING REMARKS

The main motivation was that the groups $2^3.PSL_2(7)$ and $2^3{:}PSL_2(7)$ could be used as the broken symmetry of the mass matrices of the quarks, charged leptons and the neutrinos when they are considered as a single mass matrix. The particle physics literature is full of examples of models demonstrating that the groups $PSL_2(7)$ and its maximal subgroups $7{:}3$ and $S_4$ can describe the symmetries of the neutrino mass matrix.

Here we have studied the constructions of two groups $2^3.PSL_2(7)$ and $2^3{:}PSL_2(7)$ of order 1344 and gave the character tables including those of their maximal subgroups. The simple finite subgroup $PSL_2(7)$ of $SU(3)$ occurs either as a factor group or factor and maximal subgroup in these groups. We have also seen that the extension of the Frobenious group $7{:}3$ by the elementary Abelian group $2^3$ leading to the group $2^3{:}7{:}3$ of order 168 is a subgroup in both groups of order 1344. We have also listed the tensor products of the irreducible representations of the groups $2^3.PSL_2(7)$ and $2^3{:}PSL_2(7)$ including those of their maximal subgroups in Appendix A and gave the decompositions of the irreducible representations with respect to the irreducible representations of the relevant maximal subgroups in Appendix B.

Another possible use of these groups in particle physics may arise as follows. The operators, charge conjugation $C$, parity $P$ and time-reversal $T$ generate the elementary Abelian group $2^3 = <C, P, T>$ of order 8 when they act on the bilinear Dirac fields $\bar{\psi}\Gamma\psi$ constructed by 16 $\Gamma$ matrices. The Dirac bilinear forms constitute 7 eigenvectors of the 7 operators generated by the charge conjugation, parity and the time-reversal operators with the eigenvalues $\pm 1$ besides the scalar bilinear Dirac field which represents the eigenvector of the unit operator. The representation of the generators can be taken as the matrices $N_i (i = 1,2,...,7)$. The elementary Abelian group can be extended by its automorphism group $PSL_2(7)$ to either of the group $2^3 . PSL_2(7)$ or $2^3 : PSL_2(7)$. One can construct an effective Hamiltonian as the products of Dirac bilinear fields $\bar{\psi}\Gamma\psi$. Formal properties of parity violation, $CP$ violation and even $CPT$ violation can be explained as invariances under the subgroups of the quotient group $PSL_2(7)$. If we simply impose the $\tilde{\gamma}$ invariance we can prove that the parity violating weak interaction should take the form either $V - A$ or $V + A$. If we impose only $\tilde{\beta}$ invariance we obtain $CPT$ preserving but $CP$ violating interaction. Invariance under the generator $\tilde{\alpha}$ only leads to the $CPT$ violating terms along with a necessary violation of the Lorentz invariance.

**APPENDIX A. Tensor products**

First of all, we note that the groups $7:3, 2^3:7:3, PSL_2(7), 2^3. PSL_2(7)$ and $2^3:PSL_2(7)$ all have 3-dimensional complex representations. Moreover, the group $2^3:7:3$ has three 7-dimensional irreducible representations; one real, two complex representations. The real representation preserves the octonion algebra. Tensor products of the irreducible representations of the group $PSL_2(7)$ and its maximal subgroups are listed in the reference [2], therefore we will not reproduce them. The other tensor products are given as follows.

**Tensor products of the irreducible representations of the group $2^3:7:3$**

$1_1 \times 1_1 = 1_2$
$1_1 \times 1_2 = 1$
$1_1 \times 3_1 = 3_1$
$1_1 \times 3_2 = 3_2$
$1_1 \times 7_1 = 7_3$
$1_1 \times 7_2 = 7_1$
$1_1 \times 7_3 = 7_2$

$1_2 \times 1_2 = 1_1$
$1_2 \times 3_1 = 3_1$
$1_2 \times 3_2 = 3_2$
$1_2 \times 7_1 = 7_2$
$1_2 \times 7_2 = 7_3$
$1_2 \times 7_3 = 7_1$

$3_1 \times 3_1 = 3_2 + 3_2 + 3_1$
$3_1 \times 3_2 = 1 + 1_1 + 1_2 + 3_1 + 3_2$
$3_1 \times 7_1 = 7_1 + 7_2 + 7_3$
$3_1 \times 7_2 = 7_1 + 7_2 + 7_3$
$3_1 \times 7_3 = 7_1 + 7_2 + 7_3$

$3_2 \times 3_2 = 3_1 + 3_1 + 3_2$
$3_2 \times 7_1 = 7_1 + 7_2 + 7_3$
$3_2 \times 7_2 = 7_1 + 7_2 + 7_3$
$3_2 \times 7_3 = 7_1 + 7_2 + 7_3$

$7_1 \times 7_1 = 1 + 3_1 + 3_2 + 2(7_1) + 2(7_2) + 2(7_3)$
$7_1 \times 7_2 = 1_2 + 3_1 + 3_2 + 2(7_1) + 2(7_2) + 2(7_3)$
$7_1 \times 7_3 = 1_1 + 3_1 + 3_2 + 2(7_1) + 2(7_2) + 2(7_3)$

$7_2 \times 7_2 = 1_1 + 3_1 + 3_2 + 2(7_1) + 2(7_2) + 2(7_3)$

$7_2 \times 7_3 = 1 + 3_1 + 3_2 + 2(7_1) + 2(7_2) + 2(7_3)$
$7_3 \times 7_3 = 1_2 + 3_1 + 3_2 + 2(7_1) + 2(7_2) + 2(7_3)$

**Tensor products of the irreducible representations of the groups $2^3.PSL_2(7)$ and $2^3{:}PSL_2(7)$**

$3_1 \times 3_1 = 3_2 + 6$　　　　　　　　　　$3_2 \times 3_2 = 3_1 + 6$
$3_1 \times 3_2 = 1 + 8$　　　　　　　　　　　$3_2 \times 6 = 3_1 + 7_2 + 8$
$3_1 \times 6 = 3_2 + 7_2 + 8$　　　　　　　　$3_2 \times 7_1 = 21_2$
$3_1 \times 7_1 = 21_2$　　　　　　　　　　　$3_2 \times 7_2 = 6 + 7_2 + 8$
$3_1 \times 7_2 = 6 + 7_2 + 8$　　　　　　　　$3_2 \times 7_3 = 21_1$
$3_1 \times 7_2 = 21_1$　　　　　　　　　　　$3_2 \times 8 = 3_2 + 6 + 7_2 + 8$
$3_1 \times 8 = 3_1 + 6 + 7_2 + 8$　　　　　　$3_2 \times 14 = 21_1 + 21_2$
$3_1 \times 14 = 21_1 + 21_2$　　　　　　　　$3_2 \times 21_1 = 7_3 + 14 + 21_1 + 21_2$
$3_1 \times 21_1 = 7_3 + 14 + 21_1 + 21_2$　　$3_2 \times 21_2 = 7_1 + 14 + 21_1 + 21_2$
$3_1 \times 21_2 = 7_1 + 14 + 21_1 + 21_2$

$6 \times 6 = 1 + 2(6) + 7_2 + 2(8)$
$6 \times 7_1 = 7_1 + 14 + 21_1$
$6 \times 7_2 = 3_1 + 3_2 + 6 + 2(7_2) + 2(8)$
$6 \times 7_3 = 7_3 + 14 + 21_2$
$6 \times 8 = 3_1 + 3_2 + 2(6) + 2(7_2) + 2(8)$
$6 \times 14 = 7_1 + 7_3 + 2(14) + 21_1 + 21_2$
$6 \times 21_1 = 7_1 + 14 + 3(21_1) + 2(21_2)$
$6 \times 21_2 = 7_3 + 14 + 2(21_1) + 3(21_2)$

$7_1 \times 7_1 = 1 + 6 + 7_1 + 14 + 21_1$
$7_1 \times 7_2 = 7_3 + 21_1 + 21_2$
$7_1 \times 7_3 = 7_2 + 21_1 + 21_2$
$7_1 \times 8 = 14 + 21_1 + 21_2$
$7_1 \times 14 = 6 + 7_1 + 8 + 14 + 2(21_1) + 21_2$
$7_1 \times 21_1 = 6 + 7_1 + 7_2 + 7_3 + 8 + 2(14) + 2(21_1) + 2(21_2)$
$7_1 \times 21_2 = 3_1 + 3_2 + 7_2 + 7_3 + 8 + 14 + 2(21_1) + 3(21_2)$

$7_2 \times 7_2 = 1 + 3_1 + 3_2 + 2(6) + 2(7_2) + 2(8)$
$7_2 \times 7_3 = 7_1 + 21_1 + 21_2$
$7_2 \times 8 = 3_1 + 3_2 + 2(6) + 2(7_2) + 2(8)$
$7_2 \times 14 = 14 + 2(21_1) + 2(21_2)$
$7_2 \times 21_1 = 7_1 + 7_3 + 2(14) + 2(21_1) + 3(21_2)$
$7_2 \times 21_2 = 7_1 + 7_3 + 2(14) + 3(21_1) + 2(21_2)$

$7_3 \times 7_3 = 1 + 6 + 7_3 + 14 + 21_2$
$7_3 \times 8 = 14 + 21_1 + 21_2$
$7_3 \times 14 = 6 + 7_3 + 8 + 14 + 21_1 + 2(21_2)$
$7_3 \times 21_1 = 3_1 + 3_2 + 7_1 + 7_2 + 8 + 14 + 3(21_1) + 2(21_2)$
$7_3 \times 21_2 = 6 + 7_1 + 7_2 + 7_3 + 8 + 2(14) + 2(21_1) + 2(21_2)$

$8 \times 8 = 1 + 3_1 + 3_2 + 2(6) + 2(7_2) + 3(8)$
$8 \times 14 = 7_1 + 7_3 + 14 + 2(21_1) + 2(21_2)$
$8 \times 21_1 = 7_1 + 7_3 + 2(14) + 3(21_1) + 3(21_2)$
$8 \times 21_2 = 7_1 + 7_3 + 2(14) + 3(21_1) + 3(21_2)$

$14 \times 14 = 1 + 2(6) + 7_1 + 7_2 + 7_3 + 8 + 2(14) + 3(21_1) + 3(21_2)$
$14 \times 21_1 = 3_1 + 3_2 + 6 + 2(7_1) + 2(7_2) + 7_3 + 2(8) + 3(14) + 5(21_1)$
$\qquad + 4(21_2)$
$14 \times 21_2 = 3_1 + 3_2 + 6 + 7_1 + 2(7_2) + 2(7_3) + 2(8) + 3(14) + 4(21_1) + 5(21_2)$

$21_1 \times 21_1 = 1 + 3_1 + 3_2 + 3(6) + 2(7_1) + 2(7_2) + 3(7_3) + 3(8) + 5(14) +$
$6(21_1) + 7(21_2)$

$21_1 \times 21_2 = 3_1 + 3_2 + 2(6) + 2(7_1) + 3(7_2) + 2(7_3) + 3(8) + 4(14)$
$\qquad + 7(21_1) + 7(21_2)$

$21_2 \times 21_2 = 1 + 3_1 + 3_2 + 3(6) + 3(7_1) + 2(7_2) + 2(7_3) + 3(8) + 5(14)$
$\qquad + 7(21_1) + 6(21_2).$

**Tensor products of the irreducible representations of the groups $4.S_4:2$ and $4:S_4:2$**

$1_1 \times 1_1 = 1$
$1_1 \times 1_2 = 1_3$
$1_1 \times 1_3 = 1_2$
$1_1 \times 2_1 = 2_2$
$1_1 \times 2_2 = 2_1$
$1_1 \times 3_1 = 3_4$
$1_1 \times 3_2 = 3_3$
$1_1 \times 3_3 = 3_2$
$1_1 \times 3_4 = 3_1$
$1_1 \times 6_1 = 6_2$
$1_1 \times 6_2 = 6_1$
$1_1 \times 6_3 = 6_3$
$1_1 \times 6_4 = 6_4$

$1_2 \times 1_2 = 1$
$1_2 \times 1_3 = 1_1$
$1_2 \times 2_1 = 2_2$
$1_2 \times 2_2 = 2_1$
$1_2 \times 3_1 = 3_3$
$1_2 \times 3_2 = 3_4$
$1_2 \times 3_3 = 3_1$
$1_2 \times 3_4 = 3_2$
$1_2 \times 6_1 = 6_1$
$1_2 \times 6_2 = 6_2$
$1_2 \times 6_3 = 6_4$
$1_2 \times 6_4 = 6_3$

$1_3 \times 1_3 = 1$
$1_3 \times 2_1 = 2_1$
$1_3 \times 2_2 = 2_2$
$1_3 \times 3_1 = 3_2$
$1_3 \times 3_2 = 3_1$
$1_3 \times 3_3 = 3_4$
$1_3 \times 3_4 = 3_3$
$1_3 \times 6_1 = 6_2$
$1_3 \times 6_2 = 6_1$
$1_3 \times 6_3 = 6_4$
$1_3 \times 6_4 = 6_3$

$2_1 \times 2_1 = 1 + 1_3 + 2_2$
$2_1 \times 2_2 = 1_1 + 1_2 + 2_1$
$2_1 \times 3_1 = 3_3 + 3_4$
$2_1 \times 3_2 = 3_3 + 3_4$
$2_1 \times 3_3 = 3_1 + 3_2$
$2_1 \times 3_4 = 3_1 + 3_2$
$2_1 \times 6_1 = 6_1 + 6_2$
$2_1 \times 6_2 = 6_1 + 6_2$
$2_1 \times 6_3 = 6_3 + 6_4$
$2_1 \times 6_4 = 6_3 + 6_4$

$2_2 \times 2_2 = 1 + 1_3 + 2_2$
$2_2 \times 3_1 = 3_1 + 3_2$
$2_2 \times 3_2 = 3_1 + 3_2$
$2_2 \times 3_3 = 3_3 + 3_4$
$2_2 \times 3_4 = 3_3 + 3_4$
$2_2 \times 6_1 = 6_1 + 6_2$
$2_2 \times 6_2 = 6_1 + 6_2$
$2_2 \times 6_3 = 6_3 + 6_4$
$2_2 \times 6_4 = 6_3 + 6_4$

$3_1 \times 3_1 = 1 + 2_2 + 3_1 + 3_2$  
$3_1 \times 3_2 = 1_3 + 2_2 + 3_1 + 3_2$  
$3_1 \times 3_3 = 1_2 + 2_1 + 3_3 + 3_4$  
$3_1 \times 3_4 = 1_1 + 2_1 + 3_3 + 3_4$  
$3_1 \times 6_1 = 6_2 + 6_3 + 6_4$  
$3_1 \times 6_2 = 6_1 + 6_3 + 6_4$  
$3_1 \times 6_3 = 6_1 + 6_2 + 6_4$  
$3_1 \times 6_4 = 6_1 + 6_2 + 6_3$  

$3_2 \times 3_2 = 1 + 2_2 + 3_1 + 3_2$  
$3_2 \times 3_3 = 1_1 + 2_1 + 3_3 + 3_4$  
$3_2 \times 3_4 = 1_2 + 2_1 + 3_3 + 3_4$  
$3_2 \times 6_1 = 6_1 + 6_3 + 6_4$  
$3_2 \times 6_2 = 6_2 + 6_3 + 6_4$  
$3_2 \times 6_3 = 6_1 + 6_2 + 6_3$  
$3_2 \times 6_4 = 6_1 + 6_2 + 6_4$  

$3_3 \times 3_3 = 1 + 2_2 + 3_1 + 3_2$  
$3_3 \times 3_4 = 1_3 + 2_2 + 3_1 + 3_2$  
$3_3 \times 6_1 = 6_2 + 6_3 + 6_4$  
$3_3 \times 6_2 = 6_1 + 6_3 + 6_4$  
$3_3 \times 6_3 = 6_1 + 6_2 + 6_3$  
$3_3 \times 6_4 = 6_1 + 6_2 + 6_4$  

$3_4 \times 3_4 = 1 + 2_2 + 3_1 + 3_2$  
$3_4 \times 6_1 = 6_1 + 6_3 + 6_4$  
$3_4 \times 6_2 = 6_2 + 6_3 + 6_4$  
$3_4 \times 6_3 = 6_1 + 6_2 + 6_4$  
$3_4 \times 6_4 = 6_1 + 6_2 + 6_3$  

$6_1 \times 6_1 = 1 + 1_2 + 2_1 + 2_2 + 3_2 + 3_4 + 6_1 + 6_2 + 6_3 + 6_4$  
$6_1 \times 6_2 = 1_1 + 1_3 + 2_1 + 2_2 + 3_1 + 3_3 + 6_1 + 6_2 + 6_3 + 6_4$  
$6_1 \times 6_3 = 3_1 + 3_2 + 3_3 + 3_4 + 6_1 + 6_2 + 6_3 + 6_4$  
$6_1 \times 6_4 = 3_1 + 3_2 + 3_3 + 3_4 + 6_1 + 6_2 + 6_3 + 6_4$  

$6_2 \times 6_2 = 1 + 1_2 + + 2_1 + 2_2 + 3_2 + 3_4 + 6_1 + 6_2 + 6_3 + 6_4$  
$6_2 \times 6_3 = 3_1 + 3_2 + 3_3 + 3_4 + 6_1 + 6_2 + 6_3 + 6_4$  
$6_2 \times 6_4 = 3_1 + 3_2 + 3_3 + 3_4 + 6_1 + 6_2 + 6_3 + 6_4$  

$6_3 \times 6_3 = 1 + 1_1 + 2_1 + 2_2 + 3_2 + 3_3 + 6_1 + 6_2 + 6_3 + 6_4$  
$6_3 \times 6_4 = 1_2 + 1_3 + 2_1 + 2_2 + 3_1 + 3_4 + 6_1 + 6_2 + 6_3 + 6_4$  

$6_4 \times 6_4 = 1 + 1_1 + 2_1 + 2_2 + 3_2 + 3_3 + 6_1 + 6_2 + 6_3 + 6_4$  

**Tensor products of the irreducible representations of the groups $2^3.S_4$ and $2^3{:}S_4$**

$1_1 \times 1_1 = 1$  
$1_1 \times 2 = 2$  
$1_1 \times 3_1 = 3_3$  
$1_1 \times 3_2 = 3_4$  
$1_1 \times 3_3 = 3_1$  
$1_1 \times 3_4 = 3_2$  
$1_1 \times 3_5 = 3_6$  
$1_1 \times 3_6 = 3_5$  
$1_1 \times 4_1 = 4_2$  
$1_1 \times 4_2 = 4_1$  
$1_1 \times 6 = 6$  
$1_1 \times 8 = 8$  

$2 \times 2 = 1 + 1_1 + 2$  
$2 \times 3_1 = 3_1 + 3_3$  
$2 \times 3_2 = 3_2 + 3_4$  
$2 \times 3_3 = 3_1 + 3_3$  
$2 \times 3_4 = 3_2 + 3_4$  
$2 \times 3_5 = 3_5 + 3_6$  
$2 \times 3_6 = 3_5 + 3_6$  
$2 \times 4_1 = 8$  
$2 \times 4_2 = 8$  
$2 \times 6 = 6 + 6$  
$2 \times 8 = 4_1 + 4_2 + 8$  

$$3_1 \times 3_1 = 1 + 2 + 3_1 + 3_3$$
$$3_1 \times 3_2 = 3_6 + 6$$
$$3_1 \times 3_3 = 1_1 + 2 + 3_1 + 3_3$$
$$3_1 \times 3_4 = 3_5 + 6$$
$$3_1 \times 3_5 = 3_4 + 6$$
$$3_1 \times 3_6 = 3_2 + 6$$
$$3_1 \times 4_1 = 4_1 + 8$$
$$3_1 \times 4_2 = 4_2 + 8$$
$$3_1 \times 6 = 3_2 + 3_4 + 3_5 + 3_6 + 6$$
$$3_1 \times 8 = 4_1 + 4_2 + 2(8)$$

$$3_2 \times 3_2 = 1 + 2 + 3_2 + 3_4$$
$$3_2 \times 3_3 = 3_5 + 6$$
$$3_2 \times 3_4 = 1_1 + 2 + 3_2 + 3_4$$
$$3_2 \times 3_5 = 3_3 + 6$$
$$3_2 \times 3_6 = 3_1 + 6$$
$$3_2 \times 4_1 = 4_2 + 8$$
$$3_2 \times 4_2 = 4_1 + 8$$
$$3_2 \times 6 = 3_1 + 3_3 + 3_5 + 3_6 + 6$$
$$3_2 \times 8 = 4_1 + 4_2 + 2(8)$$

$$3_3 \times 3_3 = 1 + 2 + 3_1 + 3_3$$
$$3_3 \times 3_4 = 3_6 + 6$$
$$3_3 \times 3_5 = 3_2 + 6$$
$$3_3 \times 3_6 = 3_4 + 6$$
$$3_3 \times 4_1 = 4_2 + 8$$
$$3_3 \times 4_2 = 4_1 + 8$$
$$3_3 \times 6 = 3_2 + 3_4 + 3_5 + 3_6 + 6$$
$$3_3 \times 8 = 4_1 + 4_2 + 2(8)$$

$$3_4 \times 3_4 = 1 + 2 + 3_2 + 3_4$$
$$3_4 \times 3_5 = 3_1 + 6$$
$$3_4 \times 3_6 = 3_3 + 6$$
$$3_4 \times 4_1 = 4_1 + 8$$
$$3_4 \times 4_2 = 4_2 + 8$$
$$3_4 \times 6 = 3_1 + 3_3 + 3_5 + 3_6 + 6$$
$$3_4 \times 8 = 4_1 + 4_2 + 2(8)$$

$$3_5 \times 3_5 = 1 + 2 + 3_5 + 3_6$$
$$3_5 \times 3_6 = 1_1 + 2 + 3_5 + 3_6$$
$$3_5 \times 4_1 = 4_1 + 8$$
$$3_5 \times 4_2 = 4_2 + 8$$
$$3_5 \times 6 = 3_1 + 3_2 + 3_3 + 3_4 + 6$$
$$3_5 \times 8 = 4_1 + 4_2 + 2(8)$$

$$3_6 \times 3_6 = 1 + 2 + 3_5 + 3_6$$
$$3_6 \times 4_1 = 4_2 + 8$$
$$3_6 \times 4_2 = 4_1 + 8$$
$$3_6 \times 6 = 3_1 + 3_2 + 3_3 + 3_4 + 6$$
$$3_6 \times 8 = 4_1 + 4_2 + 2(8)$$

$$4_1 \times 4_1 = 1 + 3_1 + 3_4 + 3_5 + 6$$
$$4_1 \times 4_2 = 1_1 + 3_2 + 3_3 + 3_6 + 6$$
$$4_1 \times 6 = 4_1 + 4_2 + 2(8)$$
$$4_1 \times 8 = 2 + 3_1 + 3_2 + 3_3 + 3_4 + 3_5 + 3_6 + 2(6)$$

$$4_2 \times 4_2 = 1 + 3_1 + 3_4 + 3_5 + 6$$
$$4_2 \times 6 = 4_1 + 4_2 + 2(8)$$
$$4_2 \times 8 = 2 + 3_1 + 3_2 + 3_3 + 3_4 + 3_5 + 3_6 + 2(6)$$

$$6 \times 6 = 1 + 1_1 + 2(2) + 3_1 + 3_2 + 3_3 + 3_4 + 3_5 + 3_6 + 2(6)$$
$$6 \times 8 = 2(4_1) + 2(4_2) + 4(8)$$

$$8 \times 8 = 1 + 1_1 + 2 + 2(3_1) + 2(3_2) + 2(3_3) + 2(3_4) + 2(3_5) + 2(3_6) + 4(6)$$

**APPENDIX B. Decomposition of the irreducible representations in terms of the irreducible representations of the maximal subgroups**

**Decompositions of the irreducible representations of the groups $2^3 . PSL_2(7)$ and $2^3 : PSL_2(7)$ under the maximal subgroup $2^3 : 7 : 3$**

| Irreducible representations of $2^3 . PSL_2(7)$ and $2^3 : PSL_2(7)$ | Irreducible representations of $2^3 : 7 : 3$ |
|---|---|
| 1 | 1 |
| $3_1$ | $3_1$ |
| $3_2$ | $3_2$ |
| 6 | $3_1 + 3_2$ |
| $7_1$ | $7_1$ |
| $7_2$ | $1 + 3_1 + 3_2$ |
| $7_3$ | $7_1$ |
| 8 | $1_1 + 1_2 + 3_1 + 3_2$ |
| 14 | $7_2 + 7_3$ |
| $21_1$ | $7_1 + 7_2 + 7_3$ |
| $21_2$ | $7_1 + 7_2 + 7_3$ |

**Decompositions of the irreducible representations of the group $2^3 : PSL_2(7)$ under the maximal subgroup $PSL_2(7)$**

| Irreducible representations of $2^3 : PSL_2(7)$ | Irreducible representations of $PSL_2(7)$ |
|---|---|
| 1 | 1 |
| $3_1$ | $3_1$ |
| $3_2 = \overline{3_1}$ | $3_2 = \overline{3_1}$ |
| 6 | 6 |
| $7_1$ | 7 |
| $7_2$ | 7 |
| $7_3$ | $1 + 6$ |
| 8 | 8 |
| 14 | $6 + 8$ |
| $21_1$ | $3_1 + 3_2 + 7 + 8$ |
| $21_2$ | $6 + 7 + 8$ |

**Decompositions of the irreducible representations of the group $2^3:7:3$ under the maximal subgroup $7:3$**

| Irreducible representations of $2^3:7:3$ | Irreducible representations of $7:3$ |
|---|---|
| 1 | 1 |
| $1_1$ | $1_1$ |
| $1_2 = \overline{1_1}$ | $1_2 = \overline{1_1}$ |
| $3_1$ | $3_1$ |
| $3_2 = \overline{3_1}$ | $3_2 = \overline{3_1}$ |
| $7_1$ | $1 + 3_1 + 3_2$ |
| $7_2$ | $1_2 + 3_1 + 3_2$ |
| $7_3$ | $1_1 + 3_1 + 3_2$ |

**Decompositions of the irreducible representations of the group $PSL_2(7)$ under the maximal subgroup $7:3$**

| Irreducible representations of $PSL_2(7)$ | Irreducible representations of $7:3$ |
|---|---|
| 1 | 1 |
| $3_1$ | $3_1$ |
| $3_2 = \overline{3_1}$ | $3_2 = \overline{3_1}$ |
| 6 | $3_1 + 3_2$ |
| 7 | $1 + 3_1 + 3_2$ |
| 8 | $1_1 + 1_2 + 3_1 + 3_2$ |